\documentclass[12pt]{article}

\usepackage{varioref}
\usepackage{hyperref}
\usepackage{multirow}
\usepackage{enumerate}
\usepackage{amsmath, amsfonts, amstext, amssymb, amsthm}
\usepackage[T1]{fontenc} 
\usepackage[dvips,pdftex]{graphicx}
\newtheorem{theo}{\textsc{Theorem}}[section]

\newtheorem{prop}{\textsc{Proposition}}[section]

\newtheorem{lemme}{\textsc{Lemme}}[section]

\newenvironment{preuve}[1]
{\noindent{\textbf{Proof #1:}}}{\hfill\textbf{$\square$}\newline }

\newcommand{\cE}{{\cal E}}
\newcommand{\cF}{{\cal F}}

\newcommand{\bu}{{\bf u}}

\newcommand{\by}{{\bf y}}
\newcommand{\bt}{{\bf t}}
\newfont{\msbm}{msbm10 scaled\magstep1}
\newfont{\msbms}{msbm7 scaled\magstep1} 

\newcommand{\bbC}{\mbox{$\mbox{\msbm C}$}}
\newcommand{\bbD}{\mbox{$\mbox{\msbm D}$}}

\newcommand{\bbN}{\mbox{$\mbox{\msbm N}$}}

\newcommand{\bbP}{\mbox{$\mbox{\msbm P}$}}

\newcommand{\bbR}{\mbox{$\mbox{\msbm R}$}}

\newcommand{\bbsR}{\mbox{$\mbox{\msbms R}$}}

\begin{document}

\title{Extremes of independent stochastic processes: a point process approach}
\author{Cl\'ement Dombry\footnote{Universit\'e de
Poitiers, Laboratoire de Mathématiques et Applications, UMR CNRS 6086, T\'el\'eport 2, BP 30179, F-86962 Futuroscope-Chasseneuil cedex,
France.  Email: Clement.Dombry@math.univ-poitiers.fr}\ \ and Fr\'ed\'eric Eyi-Minko\footnote{Universit\'e de
Poitiers, Laboratoire de Mathématiques et Applications, UMR CNRS 6086, T\'el\'eport 2, BP 30179, F-86962 Futuroscope-Chasseneuil cedex,
France.  Email: Frederic.Eyi.minko@math.univ-poitiers.fr}}
\date{}
\maketitle

\begin{abstract}
For each $n\geq 1$, let  $ \{ X_{in}, \quad i \geqslant 1 \} $ be  independent copies of a nonnegative continuous stochastic process $X_{n}=(X_n(t))_{t\in T}$ indexed by a compact metric space $T$. We are interested in the  process of partial maxima 
\[
\tilde M_n(u,t) =\max \{ X_{in}(t), 1 \leqslant i\leqslant [nu]  \},\quad u\geq 0,\ t\in T.
\]
where the brackets $[\,\cdot\,]$ denote the integer part. Under a regular variation condition on  the sequence of processes $X_n$, we prove that the partial maxima process $\tilde M_n$ weakly converges to a superextremal process $\tilde M$ as $n\to\infty$. We use a point process approach based on the convergence of empirical measures. Properties of the limit process are investigated: we characterize its finite-dimensional distributions, prove that it satisfies an homogeneous Markov property, and show in some cases that it is max-stable and self-similar. Convergence of further order statistics is also considered. We illustrate our results on the class of log-normal  processes in connection with some recent results on the extremes of Gaussian processes established by Kabluchko.
\end{abstract}
\ \ \ \\
{\bf Key words:} extreme value theory; partial maxima process; superextremal process;  functional regular variations; weak convergence. \\
{\bf AMS Subject classification. Primary:} 60G70\   {\bf Secondary:} 60F17
\\

\section{Introduction}
A classical problem in extreme value theory (EVT) is to determine the asymptotic behavior of the maximum of independent and identically distributed (i.i.d.) random variables $(Z_i)_{i\geq 1}$. What are the assumptions  that ensure the weak convergence of the rescaled maximum
\[
\max_{1\leq i\leq n}\frac{ Z_i-b_n}{a_n},\quad a_n>0,\ b_n\in\bbR,
\]
and what are the possible limit distributions ? These questions were at the basis of the development of EVT
and  found its answer in the Theorem by  Fisher and Tippett \cite{FT28}  characterizing all the max-stable distributions and in the description of their domain of attraction by Gnedenko \cite{G43} and de Haan \cite{dH70}. Another stimulating point of view developped by Lamperti \cite{L64} is to introduce a time variable  and consider the  asymptotic behavior of the partial maxima
\[
\max_{1\leq i\leq [nu]} \frac{ Z_i-b_n}{a_n},\quad u\geq 0.
\]
The corresponding limit process is known as an extremal process: it is a pure jump Markov process wich is also max-stable,  see Resnick \cite{R75}. 
 
Since then, EVT has known many developments. Among several other directions, the extension of the theory to multivariate and spatial settings is particularly important, as well as the statistical issues raised by the applications on real data sets. For excellent reviews of such developments, the reader is invited to refer to the monographies by Resnick \cite{L7}, de Haan and Fereira \cite{dHF06} or Beirlant {\it et al} \cite{B04} and the references therein.

Our purpose here is to focus on the functional framework and investigate the asymptotic behavior of the partial maxima processes based on a doubly infinite array of independent  random processes. Let $T$ be a compact metric space and, for each $n\geq 1$, let  $ \{ X_{in},\ i \geqslant 1 \} $ be independent copies of a sample continuous stochastic process $(X_n(t))_{t\in T}$. Without loss of generality, we will always suppose that  $X_n$ is non-negative (otherwise  consider  $X'_n(t)=e^{X_n(t)}$) and we denote by $\bbC^+=\bbC^+(T)$ be the set of non-negative continuous functions on $T$.
We are mainly interested in the process of pointwise maxima 
\begin{equation}\label{eq:defM_n}
M_n(t) =\max \{ X_{in}(t),\ 1 \leqslant i\leqslant n  \},\quad t\in T,
\end{equation}
and the process of partial maxima
\begin{equation}\label{eq:defM_ntilde}
\tilde M_n(u,t)=\max \{ X_{in}(t),\ 1 \leqslant i\leqslant [nu]  \},\quad u\geq 0,\ t\in T.
\end{equation}
We use the convention that the maximum over an empty set is equal to $0$, so $\tilde M_n(0,t)\equiv 0$. Clearly, we also have  $\tilde M_n(1,t)=M_n(t)$.
In this framework, the parameter $t\in T$ is thought as a space parameter and  $u\in [0,+\infty)$  as a time variable. 

Our approach relies on the  convergence of the following empirical measures
\begin{equation}\label{eq:defbeta_n}
\beta_{n} =  \sum_{ i = 1 }^{ n } \delta_{X_{in}} \quad \mbox{ and }\quad  \tilde \beta_{n} =  \sum_{ i \geq 1 } \delta_{ ( X_{in},i/n)}
\end{equation}
on $\bbC^+$ and $\bbC^+\times [0,+\infty)$ respectively. The maxima processes  $M_n$ and  $\tilde M_n$ can  be written as  functionals of the empirical measure  $\beta_n$ and $\tilde \beta_n$ respectively. We will prove the continuity of  the underlying functionals and then use the continuous mapping Theorem to deduce  convergence of the maxima processes  from  convergence of empirical measures.

Connections between EVT and point processes are well known. In the case when the underlying state space is locally compact, the following result holds  (see Proposition 3.21 in\cite{L7}): if $n\bbP[X_{n}\in\, \cdot\, ]$ vaguely converges to some measure $\mu$, the empirical measures $\beta_n$ and $\tilde \beta_n$  converge to Poisson random measures with intensity $\mu$ and $\mu\otimes\ell$ respectively, $\ell$ being the Lebesgue measure on $[0,+\infty)$. However, in our framework the state space $\bbC^+$ is not locally compact and we need a suitable generalization of the above result. To this aim, we follow the approach by Davis and Mikosch \cite{A4} based on the notion of boundedly finite measures and $\sharp$-weak convergence detailed by Daley and Vere-Jones in \cite{L11}.

The paper is organized as follows. In section 2, we introduce the technical material needed on  boundedly finite measures,  $\sharp$-weak convergence and convergence of empirical measures. Our main result is the convergence of the partial maxima process $ \tilde  M_{n}$ which is stated and proved in section 3. Then, in section 4, we investigate some properties of the limit process $M$, known as a superextremal process. A brief extension of our results to further order statisctics is considered in section 5.  The last section is devoted to an application of our results to  the class of log-normal processes, based on a recent work by Kabluchko \cite{A7} on the extremes of Gaussian processes.


\section{Preliminaries on boundedly finite measures and point processes}

In this section, we present general results on boundedly finite measures and point processes that will be useful in the sequel. 
The reader should refer to Appendix 2.6 in Daley and  Vere-Jones \cite{L11} or to section 2 of Davis and Mikosch \cite{A4}. This has also closed connections with the theory of regular variations, see Hult and Lindskog \cite{HL05,A5}.

\subsection{Boundedly finite measures }
Let $(E,d)$ be a complete separable metric space and  $ \mathcal{E} $  the Borel $\sigma$-algebra of $E$. 
We denote by $M_{b}(E)$ the set of all finite  measures on $(E,\cE)$. A sequence of finite measures $\{ \mu_{n}, n \geqslant 1 \}$ is said to converge weakly to $\mu\in M_{b}(E)$ if and only if $\int f d \mu_{n}  \rightarrow  \int f d \mu$ for all bounded continuous  function $f$ on $E$. We  write   $ \mu_{n} \stackrel{w}{\rightarrow} \mu$ to denote weak convergence. It is well known that this notion of convergence is metrized by the Prokhorov metric 
\[
p(\mu_1,\mu_2)= \inf \{ \varepsilon > 0;\ \forall A \in \mathcal{E}, \ \mu_1(A) \leqslant \mu_2(A^{\varepsilon} ) + \varepsilon\ \mbox{ and }\  \mu_2 (A) \leqslant \mu_1 (A^{\varepsilon} ) + \varepsilon    \}
\]
where $A^{\varepsilon}=\{x\in E;\ \exists a\in A,\ d(a,x)<\varepsilon\}$ is the $\varepsilon$-neighborhood of $A$.
Endowed with this metric, $M_{b}(E)$ is a complete separable metric space.

A  measure  $\mu$ on $(E,\cE)$  is called  boundedly finite if it assigned finite measure to bounded sets, {\it i.e.} $\mu (B) < \infty$ for all bounded $B\in \cE$. Let $M_{b}^{\sharp}(E)$ be the set of all boundedly finite measures $\mu$ on $(E,\cE)$. A sequence of boundedly finite measures $\{ \mu_{n}, n \geqslant 1 \}$ is said to converge $\sharp$-weakly to $\mu\in M_{b}^\sharp(E)$ if and only if $\int f d \mu_{n}  \rightarrow  \int f d \mu$ for all bounded continuous  function $f$ on $E$ with bounded support.
There exists a metric $p^\sharp$ on $M_{b}^{\sharp}(E)$ that is compatible with this notion of $\sharp$-weak convergence and that makes 
$M_{b}^{\sharp}(E)$ a complete and separable metric space. Such a metric can be constructed as follows. Fix an origin $e_{0} \in E $  and, for $r>0$, let $\bar B_{r}  =  \{ x \in E;\  d (x,e_0) \leq r\}$ be the closed ball of center $e_0$ and radius $r$. For any $\mu_{1}, \mu_{2}\in M_{b}^{\sharp}(E)$, let $\mu^{(r)}_{1}$ and $\mu^{(r)}_{2}$ be the restriction of $\mu_{1}$ and  $\mu_{2}$ to $\bar B_r$. Note that $\mu_{1}$ and  $\mu_{2}$ are finite measures on $\bar B_r$ and denote by $p_r$ the Prokhorov metric on $M_{b} (\bar B_{r})$. Define
\[
 p^{\sharp}(\mu_1,\mu_2) = \int_{0} ^{ \infty }  e^{-r} \frac{p_{r} (\mu_1^{ (r)} , \mu_2^{ (r) }) }{ 1 +  p_{r} (\mu_1^{ (r)} , \mu_2^{ (r) }  )} dr.
\]

There are several equivalent characterizations of $\sharp$-weak convergence. Let $\{\mu_n, n \geqslant 1 \}$ and $\mu$ be boundedly finite measures on $E$. The following statements are equivalent: 
\begin{enumerate}[ i )]
\vspace{-0.4cm}
\item  $ \int f d \mu_{n}  \rightarrow  \int f d \mu$ for all $f$ bounded continuous real valued function on $E$ with bounded support;
\item $p^{\sharp}(\mu_n,\mu)\to 0$;
\item   there exist a sequence $r_k\nearrow +\infty $ such that, for all $k \geqslant 1$,   $ \mu_{n} ^{ ( r_k )}\stackrel{w}{\rightarrow}\mu ^{ ( r_k )}$;
\item  $ \mu_{n} (B) \to \mu (B)  $ for all bounded $B \in \mathcal{E}$ such that $ \mu ( \partial B ) = 0$.
\end{enumerate}
We  write   $ \mu_{n} \stackrel{w ^{\sharp} }{\rightarrow} \mu$ to denote $\sharp$-weak convergence.

A boundedly finite point measure on $E$ is a measure $\mu$ of the  form
\[
\mu = \sum_{i\in I} \delta_{x_{i}} 
\]
with  $\{ x_{i}, i\in I \}$ a finite or countable family of points in $E$ such that any bounded set $B\subseteq E$ contains at most a finite number of the $x_i$'s. The set of boundedly finite point measures on $E$ is denoted by $M_{(b ,p)}^{\sharp}(E)$. It is a closed subset of $M_{b}^{\sharp}(E)$ and hence is complete and separable when endowed with the induced metric. 

The following characterization of $\sharp$-weak convergence of boundedly finite point measures will be useful. Let $ \{\mu_{n} = \sum_{i} \delta_{x_{i}^{n}}, n \geqslant 1\}$ and $ \mu = \sum_{i} \delta_{x_{i}}$ be elements of  $ M_{ ( b ,p )  }^{\sharp}(E)$. 
Then $\mu_{n} \stackrel{w ^{\sharp} }{\rightarrow} \mu$  if and only if there exist some sequence $r_k \nearrow + \infty $ such that for all $k\geq 1$, there exist  $N \geqslant 0$ and $m \geqslant 0$ such that
\[
\forall n \geqslant N, \quad  \mu_{n}^{(r_k)} = \sum_{i=1}^{m} \delta_{x_{i}^{n,k}}\ \mbox{ and }\ \mu^{(r_k)}  = \sum_{i=1}^{m} \delta_{x_{i}^k} 
\]
for some   points $\{x_{i}^{n,k}, x_{i}^k;\ 1\leq i\leq m,\ n\geq N\}$ in $\bar B_{r_k}$ such that $ x_{i}^{n,k} \rightarrow x_{i}^k$.

\subsection{Convergence of the empirical measures}

Let $(\Omega ,\cF,\bbP)$ be a probability space. A boundedly finite point process on $E$ is a measurable mapping 
$ N:\Omega    \rightarrow M_{ ( b ,p )  }^{\sharp}(E)$. A typical example of boundedly finite point process on $E$ is a Poisson point process $\Pi_\nu$ with intensity $\nu\in  M_{ b  }^{\sharp}(E)$. We will also consider the following empirical point processes $\beta_n$ and $\tilde\beta_n$ defined by Equation \eqref{eq:defbeta_n}, where, for each $n\geq 1$, $\{X_{in},\ i\geq 1\}$ are independent copies of an $E$-valued random variable $X_n$. The random variable $\beta_n$  is a  finite point process on $E$, while  $\tilde \beta_n$ is a boundedly finite point process on $\tilde E=E\times [0,+\infty)$ endowed with the metric
\[
 \tilde d((x_1,u_1),(x_2,u_2))=d(x_1,x_2)+|u_2-u_1|, \quad x_1,x_2\in E,\ u_1,u_2\in[0,+\infty).
\]
The following Proposition will play a key role in the sequel.

\phantomsection \begin{prop}\label{prop:pp}  
The following statements are equivalent, where the limits are taken as $n\to\infty$:
\begin{enumerate}[ i )]
\vspace{-0.4cm}
\item $n \mathbb{ P } [  X_{ n  }\in\, \cdot\, ] \stackrel{ w^{\sharp} }{\rightarrow} \nu  $;
\item $ \beta_{n}\Rightarrow \Pi_{\nu}$ with  $\Pi_\nu$ a Poisson point process on $E$ with intensity  $\mu$ and $\Rightarrow$ standing for weak convergence in $M_{ ( b ,p )  }^{\sharp}(E)$.
\item $ \tilde\beta_{n}\Rightarrow \tilde\Pi_{\nu}$ with  $\tilde\Pi_\nu$   Poisson point process on $E\times [0,+\infty)$ with intensity $\mu\otimes \ell$,  $\ell$  being the Lebesgue measure on $[0,+\infty)$ and $\Rightarrow$ standing for weak convergence in $M_{ ( b ,p )  }^{\sharp}(E\times [0,+\infty))$.
\end{enumerate}
\end{prop}
The proof of this Proposition follows by an adapation of Proposition 3.21 in Resnick \cite{L7} which states a similar result in the case of a locally compact state space and in terms of vague convergence. As noticed by  Davis and Mikosch \cite{A4}, the proof remains valid for a complete separable metric space $E$ if we change vague convergence by ${\sharp}$-weak convergence. See also Theorem 4.3 in Davydov, Molchanov and Zuyev \cite{A6}. 

In the terminology of Hult and Lindskog \cite{HL05,A5}, the condition i) in Proposition \ref{prop:pp} means that the sequence $X_n$ is {\it regularly varying}. A particularly important case is when $E$ is endowed with a strucure of cone, {\it i.e.} a multiplication by positive scalars. 
If there exists a sequence $(a_n)_{n\geq 1}$ of positive reals such that 
\[
n \mathbb{ P } [  a_n^{-1} X_{1  }\in \,\cdot\, ] \stackrel{ w^{\sharp} }{\rightarrow} \nu
\]
for some nonzero $\nu\in M_b^\sharp(E)$, then the random variable $X_1$ is said to be {\it regularly varying}. Under some technical assumptions on the structure of the convex cone $E$ ({\it e.g.} continuity properties of the multiplication), the limit measure $\nu$ is proved to be homogeneous of order $-\alpha<0$, {\it i.e.}
\[
 \nu(\lambda A)=\lambda^{-\alpha} \nu(A),\quad \lambda >0, A\in\cE.
\]
Furthermore, the sequence $a_n$ is regularly varying with index $1/\alpha$. In this framework, if $\{X_{i},\ i\geq 1\}$ is an i.i.d. sequence of regularly varying random variables and $X_{in}=a_n^{-1}X_{i}$, then condition i) is equivalent to  
\[
n \mathbb{ P } [  a_n^{-1} X_{1  }\in \,\cdot\, ] \stackrel{ w^{\sharp} }{\rightarrow} \nu
\]
where the limit measure $\nu$ is homogeneous of order $-\alpha<0$ and 
\[
M_n(u,t)=a_n^{-1} \max_{1\leq i\leq [nu]} X_i(t),\quad u\geq 0,\ t\in T.
\]

\section{Extremes of independent stochastic processes}

We go back to our original problem where the $X_{in}$'s are non negative continuous processes on $T$. We endow the set $\bbC^+=\bbC^+(T)$ of non negative continuous functions on $T$ with the uniform norm  $\|x\|=\sup_{t\in T} |x(t)|$, $x\in\bbC^+$. It turns out that, when working with maxima, the uniform metric is not adaptated, mainly because sets such as $\{x\in\bbC^+; \|x\|\geq \varepsilon\}$ are not bounded for this metric. For this reason, we introduce 
$ \overline{ \mathbb C }_{0}^{+} = ( 0 , + \infty ] \times \mathbb{S}_{\mathbb{C}^{+}}  $ 
where $\mathbb{S}_{\mathbb{C}^{+}}= \{ x \in \mathbb C^{+} : \| x \|_{\infty} = 1 \}$ is the unit sphere.
We define the metric
\[
d ( ( r_{1}, s_{1}) ,( r_{2}, s_{2})   ) = | 1/ r_{1}  - 1/r_{2}  | + \| s_{1} - s_{2} \| , \quad    ( r_{1}, s_{1})  ,  (  r_{2}  , s_{2}   ) \in  \overline{ \mathbb C}_{0}^{+}. 
\]
The metric space $(\overline{ \mathbb C }_{0}^{+},d)$ is complete and separable. Define also $\bbC_0^+=\bbC^+\setminus\{0\}$ and consider the ``polar decomposition'':
\[
T : \left\{\begin{array}{ccc} \mathbb C_{0}^{+}&\to&  ( 0 , \infty ) \times \mathbb{S}_{\mathbb{C}^{+}}  \\ x& \mapsto  &( \| x \| , x / \| x \| ) \end{array}\right..
\]
The mapping $T$ is an homeomorphism and we identify in the sequel $\mathbb C_{0}^{+}$ and $( 0 , + \infty ) \times \mathbb{S}_{\mathbb{C}^{+}}$. In this metric, a subset $B$ of $\bbC^+_0$ is bounded if and only if it is bounded away from zero, {\it i.e.} included in a set of the form $\{x\in\bbC^+; \|x\|\geq \varepsilon\}$ for some $\varepsilon>0$.

\subsection{Spatial maximum process}

For the sake of clarity, we present first our results on the spatial maximum process $M_n$ defined by Equation \eqref{eq:defM_n}. 
\begin{theo}\label{theo1}
Assume that $n \mathbb{ P } [   X_{n}  \in \,\cdot\, ] \stackrel{ w^{\sharp} }{\rightarrow} \nu $ in $M_b^\sharp(\overline\bbC_0^+)$ and that $\nu(\overline\bbC_0^+\setminus\bbC_0^+)=0$. Then, the process $M_n$ weakly converges   in $\bbC(T)$ as $n\to\infty$ to  the process $M$ defined by
\begin{equation}\label{eq:theo1}
M(t)=\sup_{i\geq 1}Y_{i}(t),\quad t\in T,
\end{equation}
where $\sum_{i\geq 1} \delta_{Y_{i}}$ is a Poisson point process on $\bbC_0^+$ with intensity $\nu$.
\end{theo}

The proof of Theorem \ref{theo1} relies on the following Lemma. With a slight abuse of notation, we define
\[
{M}_{ ( b ,p )  }^{\sharp}(\bbC_0^+) = \{ \mu \in {M}_{ ( b ,p )  }^{\sharp}(\overline{\bbC}_0^+) : \mu ( \overline{\bbC}_0^+\setminus\bbC_0^+)= 0 \}.
\]
Note ${M}_{ ( b ,p )  }^{\sharp}(\bbC_0^+)$ is an open subset of ${M}_{ ( b ,p )  }^{\sharp}(\overline{\bbC}_0^+)$.

\begin{lemme}\label{lemme1}
The mapping $\theta: {M}_{ ( b ,p )  }^{\sharp}(\bbC_0^+)\to\bbC^+$ defined by $\theta(0)\equiv 0$ and 
\[
 \theta\Big(\sum \limits_{  i \in I } \delta_{ ( r_{i} , s_{i} ) }\Big)=\Big( t \mapsto \sup\{ r_{i} s_{i} (t);\  i \in I\} \Big)
\]
is well-defined and continuous.
\end{lemme}

\begin{preuve}{of Lemma \ref{lemme1}}\\
\begin{itemize}
\vspace{-0.4cm}
 \item First, we show that the mapping $\theta$ is well defined.  This is not obvious since the pointwise supremum of a countable family of functions is not necessarily finite nor continuous.\\
 For $\varepsilon>0$,  let $\mu^{ ( \varepsilon ) }$ be the restriction of $\mu$ on the set $\bar B^{\varepsilon}=[\varepsilon , \infty ] \times \mathbb{S}_{\mathbb{C}} $. Since  $\bar B^{  \varepsilon  }$ is bounded for the metric $d$, the point measure $ \mu^{ ( \varepsilon ) } $ has only a finite number of atoms and therefore $ \theta( \mu^{ ( \varepsilon ) } )$ is the maximum of a finite number of non negative continuous functions and is hence a non negative continuous function. 

Furthermore, for all $t \in T$, $|\theta(\mu^{ ( \varepsilon ) })(t)-\theta(\mu)(t)| \leqslant \varepsilon$. So $\theta(\mu)$ is the uniform limit as $\varepsilon\to 0$ of the continuous functions $\theta _{\varepsilon }  ( \mu)$, and hence $\theta(\mu)\in\bbC^+$ and satisfies 
\begin{equation}\label{eq:estim1}
\| \theta( \mu^{ ( \varepsilon ) }) - \theta  ( \mu) \| \leqslant \varepsilon.
\end{equation}
\item Second, we show that the mapping $\theta$ is continuous.\\
Let $ \{ \mu_{n}= \sum_{i \in I_{n}} \delta_{x_{i}^{n}} , \quad n \geqslant 1 \} $ be a sequence of measures in  ${M}_{ ( b ,p )  }^{\sharp}(\bbC_0^+)$    converging to $\mu = \sum_{i \in I} \delta_{x_{i}} $. For all $\varepsilon > 0$, there exists $\varepsilon'<\varepsilon$,  $N \geqslant 0$ and  $m \geqslant 0$ such that
\[
 \forall n \geqslant N, \quad  \mu_{n}^{( \varepsilon' )} = \sum_{i=1}^{m} \delta_{x_{i}^{n\varepsilon'}}\ \mbox{ and }\ \mu^{( \varepsilon' )}  = \sum_{i=1}^{m} \delta_{x_{i}^{\varepsilon'}}
\]
with $x_{i}^{n\varepsilon'}, x_{i}^{\varepsilon'} \in \bar B^{  \varepsilon'  } $ and  $ x_{i}^{n\varepsilon'} \rightarrow x_{i}^{\varepsilon'}$ as $n \rightarrow \infty$. Clearly, this implies that $\theta(\mu_{n}^{( \varepsilon')})$ converges in $\bbC^+$ to $\theta(\mu^{( \varepsilon')})$ as $n \rightarrow \infty$. 
Hence there exists an integer $n_0$ such that    $\|\theta ( \mu_{n}^{ ( \varepsilon' ) }) - \theta  ( \mu^{ ( \varepsilon ') }) \| \leq \varepsilon$ for all $n\geq n_0$.\\
Then, thanks to Equation \eqref{eq:estim1}, we obtain for $n\geq n_0$
\[
\|\theta  ( \mu_{n})-\theta  ( \mu)\|\leq  
 \| \theta  ( \mu_{n}) -\theta  ( \mu_{n}^{ ( \varepsilon' ) })  \| + \|\theta  ( \mu_{n}^{ ( \varepsilon' ) })- \theta  ( \mu^{ ( \varepsilon ') }) \| + \| \theta( \mu^{ ( \varepsilon' ) }) - \theta  ( \mu) )\|\leq 3\varepsilon. 
\]
This proves that $\theta(\mu_n)\to \theta(\mu)$  as $n\to\infty$ and that the mapping $\theta$ is continuous.
\end{itemize}
\end{preuve}

\begin{preuve}{of Theorem \ref{theo1}}\\
Under the assumption  $n \mathbb{ P } [   X_{n}\in \,\cdot\, ] \stackrel{ w^{\sharp} }{\rightarrow} \nu $, we know from Proposition \ref{prop:pp} that the empirical measure $\beta_{n}$ defined by Equation \eqref{eq:defbeta_n} 
weakly converges in $M_{(b,p)}(\overline{\bbC}_0^+)$ to a Poisson point process $\Pi_\nu$ with intensity $\nu$. The assumption $\nu(\overline{\bbC}_0^+\setminus \bbC_0^+)$ ensures that $\Pi_\nu$ lies almost surely in $M_{(b,p)}(\bbC_0^+)$. Note that $M_n=\theta(\beta_n)$, and according to Lemma \ref{lemme1}, the mapping $\theta$ is continuous. So the continuous mapping Theorem (see {\it e.g.} Theorem 5.1 in Billingsley \cite{B68})  entails $\theta ( \beta_{n} ) \Rightarrow \theta ( \Pi_{\nu} )$ which is equivalent to $M_n \Rightarrow M$.
\end{preuve}

\subsection{Spatio-temporal maximum process}

We consider now convergence of the space-time process $\tilde M_n$ defined by Equation \eqref{eq:defM_ntilde}. 

For fixed $u\geq 0$, the space process $t\mapsto \tilde M_n(u,t)$ is sample continuous and non negative, i.e. a random elements of $\bbC^+$. Furthermore, the time process $u\mapsto \tilde M_n(u,\cdot)$ can be seen as a $\bbC^+$-valued càd-làg process on $[0,+\infty)$; it is indeed constant on intervals of the form $[k/n,k/n+1/n)$, $k\in\bbN$. Hence, we will consider the process $\tilde M_n$ as a random element of the Skohorod space $\bbD([0,+\infty),\bbC^+)$ endowed with the $J_1$-topology (see for example Ethier and Kurtz \cite{EK86} for the definition and properties of Skohorod space). 

\begin{theo}\label{theo2}
Assume that $n \mathbb{ P } [   X_{n}  \in \,\cdot\, ] \stackrel{ w^{\sharp} }{\rightarrow} \nu $ in $M_b^\sharp(\overline\bbC_0^+)$ and that $\nu(\overline\bbC_0^+\setminus\bbC_0^+)=0$. Then, the process $\tilde M_n$ weakly converges   in $\bbD([0,+\infty),\bbC^+)$ as $n\to\infty$ to the superextremal process $\tilde M$ defined by
\begin{equation}\label{eq:theo2}
\tilde M(u,t)=\sup\{Y_{i}(t)\mathbf 1_{ [U_i,+\infty)}(u) ;\ i\geq 1\}, \quad u>0,\ t\in T,
\end{equation}
where $\sum_{i\geq 1} \delta_{(Y_{i},U_i)}$ is a Poisson point process on $\bbC_0^+\times [0,+\infty)$ with intensity $\nu\otimes\ell$.
\end{theo}

For the proof, we will need the following analoguous of Lemma \ref{lemme1} in the space-time framework. Let $M_{b,p}^\sharp(\bbC_0^+\times [0,+\infty))$ be the subset of measures $\mu\in M_{(b,p)}^\sharp(\overline{\bbC}_0^+\times [0,+\infty))$
such that $\mu((\overline{\bbC}_0^+\setminus \bbC_0^+)\times [0,+\infty))=0$ and define
\[
\tilde C=\Big\{\mu\in M_{(b,p)}^\sharp(\bbC_0^+\times [0,+\infty));\ \mu(\overline{\bbC}_0^+\times\{u\})\leq 1\quad \mathrm{for\ all\ } u\geq 0 \Big\}.
\]
In other words, a measure $\mu=\sum_{i\in I}\delta_{ (  r_{i}, s_{i}  , u_{i} ) }$ belongs to $\tilde C$ if and only if  $r_i<+\infty$ for all $i\in I$ and the $u_i$'s are pairwise distinct.

\begin{lemme} \label{lemme2}
The mapping $\tilde\theta:{M}_{ ( b ,p )  }^{\sharp}(\bbC_0^+\times[0,+\infty))\to \bbD([0,+\infty),\bbC^+)$ defined by $\tilde\theta(0)\equiv 0$ and 
\[
\tilde \theta\Big(\sum \limits_{  i \in I } \delta_{ ( r_{i} , s_{i},u_i ) }\Big)=\Big( (u,t) \mapsto \sup\{   r_{i} s_{i} (t)1_{[ u_{i} , + \infty )} (u);\ i\in I\} \Big)
\]
is well defined and continuous on $\tilde C$.
\end{lemme}

\begin{preuve}{of Lemma \ref{lemme2}}\\
The proof is similar to the proof of Lemma \ref{lemme1} and we give only the main lines.
\begin{itemize}
\vspace{-0.4cm}
\item First we show that $\tilde\theta$ is well defined. Recall that $\overline{\bbC}_0^+\times[0,+\infty)$ is endowed with the metric \[
\tilde d((r_1,s_1,u_1),(r_2,s_2,u_2))=d((r_1,s_1),(r_2,s_2))+|u_2-u_1|.
\]
 Let $\mu\in {M}_{ ( b ,p )  }^{\sharp}(\bbC_0^+\times[0,+\infty))$. For $\varepsilon>0$ and $M>0$, let $\mu^{(\varepsilon,M)}$ be its restriction to $\bar B^{\varepsilon,M}=[\varepsilon, +\infty ] \times \mathbb{S}_{\mathbb{C}^{+}} \times [ 0, M]$. Since $\bar B^{\varepsilon,M}$ is bounded, $\mu^{(\varepsilon,M)}$ has only a finite number of atoms and we easily check that $\tilde \theta(\mu^{(\varepsilon,M)})$ belongs to $\bbD([0,+\infty),\bbC^+)$. Furthermore,  for $u<M$ and $t\in T$,
\[
|\tilde\theta(\mu^{(\varepsilon,M)})(u,t)-\tilde\theta(\mu)(u,t)|\leq \varepsilon 
\]
so that $\tilde\theta(\mu^{(\varepsilon,M)})$ converges uniformly on $[0,M]\times T$ to $\tilde\theta(\mu)$ as $\varepsilon\to 0$. The constant $M$ being arbitrary, this implies that $\tilde\theta(\mu)\in\bbD([0,+\infty),\bbC^+)$ and that the application $\tilde\theta$ is well defined.
It also holds that
\begin{equation}\label{eq:estim2}
 \sup_{(u,t)\in[0,M]\times T} |\tilde\theta(\mu^{(\varepsilon,M)})(u,t)-\tilde\theta(\mu)(u,t)|\leq \varepsilon.
\end{equation}
\item Next, we show that $\tilde\theta$ is continuous on $\tilde C$.\\
 Let $ \{ \mu_{n}= \sum_{i \in I_{n}} \delta_{x_{i}^{n}}; \ n \geqslant 1 \} $ be a sequence of $M_{ ( b ,p )  }^{\sharp}(\bbC_0^+\times[0,+\infty))$ converging to $\mu = \sum_{i \in I} \delta_{x_{i}} \in \tilde C$.
For all $\varepsilon>0$ and $M>0$, there exist $\varepsilon' < \varepsilon$, $M'>M$, $N \geq 1$ and some $m\geq 1$, such that
\[
 \forall n \geqslant N, \quad  \mu_{n}^{( \varepsilon',M')} = \sum_{i=1}^{m} \delta_{x_{i}^{n}}\ \mbox{ and }\ \mu^{( \varepsilon',M')}  = \sum_{i=1}^{m} \delta_{x_{i}}
\]
with  $x_{i}^{n} = (r_{i}^n , s_{i}^n  , u_{i}^n )$, $x_{i} = (r_{i} , s_{i}  , u_{i} )\in  \bar B^{  \varepsilon',M'  }$ and $ x_{i}^{n} \to x_{i}$ as $n \rightarrow \infty$.
The condition $\mu\in \tilde C$ ensures that the $u_i$'s are pairwise distinct and we can suppose without loss of generality that $u_1<\cdots<u_m$. For $n$ large enough, we will also have
$u_1^n<\cdots<u_m^n$. Define $\delta_{M'}$ the metric associated with the $J_1$-topology on $\bbD([0,M'],\bbC^+)$  by
\[
 \delta_{M'}(x,y)=\inf_\lambda \sup_{(u,t)\in[0,M']\times T} |x(\lambda u,t)-x(u,t)|
\]
where the  infinimum is taken over the set of non-decreasing  homeomorphisms $\lambda$ of $[0,M']$. Since for large $n$ the $u_i$'s and the $u_i^n$'s are in the same relative order, there exists a non-decreasing homeomorphism $\lambda_{M'}^n$ of $[0,M']$ such that $\lambda_{M'}^n(u_i^n)=u_i$ for all $1\leq i\leq m$. We then have 
$$
\begin{tabular}{rcl}
\multicolumn{3}{l}{$\delta_{M'}(\tilde\theta(\mu_n^{(\varepsilon',M')}),\tilde\theta(\mu^{(\varepsilon',M')}))$} \\
 &$ \leq $ & $\delta_{M'}(\tilde\theta(\mu_n)) ,\tilde\theta(\mu_n^{(\varepsilon',M')}) )+\delta_{M'}(\tilde\theta(\mu_n^{(\varepsilon',M')}) ,\tilde\theta(\mu^{(\varepsilon',M')}) )+\delta_{M'}(\tilde\theta(\mu^{(\varepsilon',M')}) ,\tilde\theta(\mu) )$\\

& $ \leq$ & $\sup \limits_{(u,t)\in[0,M']\times T} |\max_{1\leq i\leq m}(r_i^ns_i^n(\cdot)\mathbf 1_{[\lambda_{M'}^n u_i^n,+\infty)}(u))-\max \limits_{1\leq i\leq m}(r_is_i(\cdot)\mathbf 1_{[u_i,+\infty)}(u))|$\\

&$=$& $\max\limits_{1\leq i\leq m} \|r_i^ns_i^n-r_is^i\|$\\
&$\to$&$ 0 \quad \mathrm{as}\ n\to\infty.$

\end{tabular}
$$
Hence, for sufficiently large $n$, $\delta_{M'}(\tilde\theta(\mu_n^{(\varepsilon',M')}) ,\tilde\theta(\mu^{(\varepsilon',M')}) )\leq \varepsilon$ and thanks to Equation \eqref{eq:estim2}, this entails
$$
\begin{tabular}{rcl}
\multicolumn{3}{l}{$\delta_{M'}(\tilde\theta(\mu_n)) ,\tilde\theta(\mu) )$} \\
 &$ \leq $ & $\delta_{M'}(\tilde\theta(\mu_n)) ,\tilde\theta(\mu_n^{(\varepsilon',M')}) )+\delta_{M'}(\tilde\theta(\mu_n^{(\varepsilon',M')}) ,\tilde\theta(\mu^{(\varepsilon',M')}) )+\delta_{M'}(\tilde\theta(\mu^{(\varepsilon',M')}) ,\tilde\theta(\mu) )$\\
 &$\leq$& 3$\varepsilon.$
\end{tabular}
$$
Since $M'$ is arbitrary large and $\varepsilon$ arbitrary small, this proves the convergence $\tilde\theta(\mu_n)$ to $\tilde\theta(\mu)$ in $\bbD([0,+\infty),\bbC^+)$.
\end{itemize}
\end{preuve}

\begin{preuve}{ of Theorem \ref{theo2} }\\
The proof is very similar to the proof of Theorem \ref{theo1}. Note that $\tilde M_n=\tilde\theta(\tilde\beta_n)$ where $\tilde\beta_n$ is defined by Equation \eqref{eq:defbeta_n}. According to Proposition \ref{prop:pp}, the empirical measure
$\tilde\beta_{n}$ weakly converges in $M_{(b,p)}(\overline{\bbC}_0^+\times [0,+\infty))$ to a Poisson point process $\tilde\Pi_\nu$ with intensity $\nu\otimes\ell$. The assumption $\nu(\overline{\bbC}_0^+\setminus \bbC_0^+)$ and the fact that the Lebesgue measure $\ell$ has no atoms ensure that $\tilde\Pi_\mu$ lies almost surely in $\tilde C$. From Lemma \ref{lemme2},  the mapping $\tilde\theta$ is continuous on $\tilde C$, so that the continuous mapping Theorem implies $\tilde\theta(\tilde\beta_n)\Rightarrow \tilde\theta(\tilde\Pi_\nu)$ which is equivalent to $\tilde M_n\Rightarrow\tilde M$. 
\end{preuve}

\section{Properties of the limit process}

In this section, we give some properties of the superextremal process $\tilde M$ defined by Equation \eqref{eq:theo2} in Theorem \ref{theo2}. First, we characterize its finite dimensional distributions. We use vectorial notations: for $l\geq 1$,  $\bt=(t_1,\ldots,t_l)\in T^l$, $\by=(y_1,\ldots,y_l)\in [0,+\infty)^l$ and $u>0$, we write $\tilde M(u,\bt)\leq \by$ if and only if $\tilde M(u,t_i)\leq y_i$ for all $i\in\{1,\ldots,l\}$.

\begin{prop}\label{prop2}
For $k \geqslant 1$, $0=u_0<u_1<\cdots<u_k$, $\bt \in T^l$ and $\by_1,\ldots,\by_k\in [0,+\infty)^l$, it holds
\[
 \bbP \Big[ \tilde M (u_{1},\bt)\leqslant \by_{1},... , \tilde M (  \bu_{k} , \bt   )\leqslant \by_{k}\Big]
= \prod_{i=1}^{k}  \exp [ -(u_j-u_{j-1})   \nu(A_j)]
\]
with  $A_j=\{f\in\bbC_0^+;\   \exists i\in 1,\ldots,l,\ f(t_{i}) \geq \min_{k\geq j} y_{k,i} \}.$
\end{prop}

\begin{preuve}{of Proposition \ref{prop2} }\\
Let $j\in\{1,\ldots,k\}$. Note that $\tilde M (u_j,\bt)\leqslant \by_j$ if and only if $\tilde\Pi_\mu$ doesn't intersect the set  
\[
B_j= \{(f,u)\in \bbC_0^+\times [0,+\infty);\ u\leq u_j \ \mathrm{and}\ \exists i\in \{1,\cdots l\}\ f(t_i)>y_{j,i}\}
\]
So, we have
\begin{eqnarray*}
\bbP \Big[ \tilde M (u_{1},\bt)\leqslant \by_{1},... , \tilde M (  \bu_{k} , \bt   )\leqslant \by_{k}\Big]
&=& \bbP \Big[ \tilde\Pi_\nu \cap (\cup_{j=1}^k  B_j)=\emptyset]\\
&=&\exp[-(\nu\otimes\ell)(\cup_{j=1}^k  B_j)].
\end{eqnarray*}
The $B_j$'s are not pairwise disjoint. To compute the measure $(\nu\otimes\ell)(\cup_{j=1}^k  B_j)$, we observe that 
\[
 \cup_{j=1}^k  B_j=\cup_{j=1}^k  ( [u_{j-1},u_j)\times A_j) \cup (\{u_k\}\times A_k)
\]
where the sets in the right hand side are pairwise  disjoint. From this, we deduce
\[
 (\nu\otimes\ell)(\cup_{j=1}^k  B_j)=\sum_{j=1}^k (u_j-u_{j-1})   \nu(A_j).
\]
This proves the Proposition.
\end{preuve}

Next we prove that the process $u\mapsto \tilde M(u,\cdot)$ is a $\bbC^+$-valued homogeneous Markov process. Let $\mathcal F_{u}$ be the $\sigma$-algebra generated by $\{\tilde M (s,t);\  s\in [0,u], t \in T\}$. The symbol $\vee$ stands for poinwise maximum.
\begin{prop}\label{prop3}
Let $u\geq 0$. The conditional distribution of $(\tilde M(u +h,\cdot))_{h\geq 0}$ given $\cF_u$ is equal to the distribution of $(\tilde M(u,\cdot)\vee\tilde M'(h,\cdot))_{h\geq 0}$ where $\tilde M'$ is an independent copy of $\tilde M$. 
\end{prop}

In some sense, this Proposition states that the process $\tilde M$ has ``independent and stationary increments'' with respect to the maximum: for $0=u_0<u_1<\cdots<u_k$, the distribution of $(\tilde M(u_i,\cdot))_{1\leq i\leq k}$ is equal to the distribution of 
\[
(\vee_{j=1}^i \tilde M^j(u_j-u_{j-1},\cdot))_{1\leq i\leq k},
\]
where $\tilde M^1,\ldots,\tilde M^k$ are i.i.d. copies of $\tilde M$. This property is similar to the property of independence and stationarity of increments (with respect to the addition) of L\'evy processes. For a fixed point  $t\in T$, the process $ \{ M ( u , t), \quad u \geq 0\}$ is known as an extremal process (see Proposition 4.7 of Resnick \cite{L7}).

\begin{preuve}{of Proposition \ref{prop3}}\\
Consider the decomposition $\tilde \Pi_\mu=\tilde\Pi_\nu^{[0,u]}\cup \tilde\Pi_\nu^{(u,\infty)}$
where
\[
\tilde\Pi_\nu^{[0,u]}=\tilde \Pi_\nu\cap (\bbC_0^+\times [0,u])\quad \mbox{and}\quad 
\tilde\Pi_\nu^{(u,\infty)}=\tilde \Pi_\nu\cap (\bbC_0^+\times (u,\infty)).
\]
By the independence properties of the Poisson point process, $\tilde\Pi_\mu^{[0,u]}$ and $\tilde\Pi_\mu^{(u,\infty)}$
are independent. Furthermore,
\begin{eqnarray*}
\tilde M(u+h,\cdot) &=&\sup\{rs(\cdot)1_{[v,+\infty)}(u+h);\ (r,s,v)\in\tilde\Pi_\nu^{[0,u]}\} \\
&&\quad \vee \ \sup\{rs(\cdot)1_{[v,+\infty)}(u+h);\ (r,s,v)\in\tilde\Pi_\nu^{(u,\infty)}\}
\end{eqnarray*}
and the two terms in the r.h.s. are independent. It is easily seen that the first term is equal to $\tilde M(u,\cdot)$, and that the invariance of the Lebesgue measure $\ell$ implies that the second term has the same distribution as $\tilde M(h,\cdot)$. This proves the Proposition.
\end{preuve}

In the particular case when the measure $\nu$ is homogeneous of order $-\alpha<0$, the process $\tilde M$ enjoys further interesting properties.
\begin{prop}\label{prop4}
Suppose $\nu$ is homogeneous with index $-\alpha < 0$. Then:
\begin{itemize}
\vspace{-0.4cm}
\item $\tilde M$ is max-stable with index $\alpha>0$, i.e. if $\tilde M^{1},\ldots, \tilde M^{n}$ are independent copies of $\tilde M$, the maximum $\vee_{i=1}^n \tilde M^i$ has the same law as $ n^{\frac{1}{\alpha}} \tilde M$. 
\item $\tilde M$ is  self-similar with index $\frac{1}{\alpha}$, i.e. for all $c>0$, the rescaled process $(\tilde M(cu,\cdot))_{u\geq 0}$ has the same distribution as $(c^{1/\alpha}\tilde M(cu,\cdot))_{u\geq 0}$
\end{itemize}
\end{prop}

\begin{preuve}{of Proposition \ref{prop4} }\\
We check that the the two processes have the same finite dimensional distributions. To this aim, we use Proposition \ref{prop2} and the notations of the proposition.
\begin{itemize}
\vspace{-0.4cm}
\item Let $\tilde M^{1},\ldots, \tilde M^{n}$ be independent copies of $\tilde M$. By Proposition \ref{prop2} and the homogeneity of $\mu$, 
\begin{eqnarray*}
 \bbP \Big[ \vee_{i=1}^n \tilde M^i (u_{j},\bt)\leqslant \by_{j},\ \forall j\in\{1,\ldots,k\}\Big]
&=& \bbP \Big[ \tilde M (u_{j},\bt)\leqslant \by_{j},\ \forall j\in\{1,\ldots,k\}\Big] ^n\\
&=& \prod_{i=1}^{k}  \exp [ -(u_j-u_{j-1})   \nu(A_j)]^n\\
&=& \prod_{i=1}^{k}  \exp [ -(u_j-u_{j-1})   \nu(n^{-1/\alpha} A_j)]\\
&=& \bbP \Big[ n^{1/\alpha}\tilde M (u_{j},\bt)\leqslant \by_{j},\ \forall j\in\{1,\ldots,k\}\Big]. 
\end{eqnarray*}
This proves the max-stability.
\item Similarly, the self-similarity is proven as follows:
\begin{eqnarray*}
 \bbP \Big[  \tilde M (cu_{j},\bt)\leqslant \by_{j},\ \forall j\in\{1,\ldots,k\}\Big]
&=& \prod_{i=1}^{k}  \exp [ -c(u_j-u_{j-1})   \nu(A_j)]\\
&=& \prod_{i=1}^{k}  \exp [ -(u_j-u_{j-1})   \nu(c^{-1/\alpha} A_j)]\\
&=& \bbP \Big[ c^{1/\alpha}\tilde M (u_{j},\bt)\leqslant \by_{j},\ \forall j\in\{1,\ldots,k\}\Big]. 
\end{eqnarray*}
\end{itemize}
\end{preuve}

\section{Further order statistics}
The point process approach for extremes is powerful: the convergence of the empirical measure $\beta_n$ entails not only the convergence of the maxima $M_n$ but also of all the order statistics, and a similar result hold for the space-time version of these processes. 

For $1\leq r\leq n$ and $t\in T$, define $M_n^r(t)$ as the $r$-th largest value among 
$X_{1n}(t),\ldots,X_{nn}(t)$. Note $M_n^1(t)= M_n(t)$ is simply the maximum. For $r>n$, we use the convention $M_n^r(t)= 0$.
The process $M_n^r=(M_n^r(t))_{t\in T}$ is refered to as the $r$-th order statistic of the sample $\{X_{1n},\ldots,X_{1n}\}$.
Similarly, for $r\geq 1$, $t\in T$ and $u\geq 0$ define $\tilde M_n^r(u,t)$ as the $r$-th largest value among $X_{1n}(t),\ldots,X_{[nu]n}(t)$ with the convention $\tilde M_n^r(u,t)=0$ if $r>[nu]$. An alternating definition in terms of the empirical measure $\tilde\beta_n$ is given by
\[
 \tilde M_n^r(u,t)=\sup\{y\geq 0;\  \tilde\beta_n(B_{u,t}^y)\geq r\} 
\]
with $B_{u,t}^y=\{(f,v)\in \bbC_0^+\times [0,+\infty); f(t)\geq y, v\leq u\}$ and the convention that the supremum over an empty set is equal to zero.

With these notations, we can strengthen Theorem \ref{theo2} as follows:

\begin{theo}\label{theo3}
Assume that $n \mathbb{ P } [   X_{n}  \in \,\cdot\, ] \stackrel{ w^{\sharp} }{\rightarrow} \nu $ in $M_b^\sharp(\overline\bbC_0^+)$ and that $\nu(\overline\bbC_0^+\setminus\bbC_0^+)=0$. Then, for each $r\geq 1$, the joint order statistics $(\tilde M_n^1,\ldots, \tilde M_n^r)$ weakly converges   in $\bbD([0,+\infty),(\bbC^{+})^r)$ as $n\to\infty$ to  $(\tilde M^1,\ldots,\tilde M^r)$ defined by
\[
\tilde M^j(u,t)=\sup\{y\geq 0;\  \tilde\Pi_\nu\cap B_{u,t}^y\geq r\} , \quad 1\leq j\leq r,\ u\geq 0,\ t\in T,
\]
with $\tilde\Pi_\nu$  a Poisson point process on $\bbC_0^+\times [0,+\infty)$ with intensity $\nu\otimes\ell$.
\end{theo}
The limit process $\tilde M^j$ corresponds to the $r$-th space-time order statistic associated to the point process $\tilde \Pi_\nu$.

\begin{preuve}{of Theorem \ref{theo3} }\\
Using once again Proposition \ref{prop:pp} and the continuous mapping Theorem, it is enough to prove the following generalization of Lemma \ref{lemme2}: the mapping 
\[
\tilde\theta^r: {M}_{ ( b ,p )  }^{\sharp}(\bbC_0^+\times[0,+\infty))\to \bbD([0,+\infty),(\bbC^+)^r) 
\]
defined by $\tilde\theta(0)\equiv 0$ and 
\[
\tilde \theta (\beta)(u,t)=(\sup\{y\geq 0;\  \tilde\beta(B_{u,t}^y)\geq j\})_{1\leq j\leq r}
\]
is well defined and continuous on $\tilde C$. The proof is very similar to the proof of Lemma \ref{lemme2} and  the details are omitted for the sake of brevity.
\end{preuve}

\section{Gaussian case}
In geostatistics, Gaussian processes are often used as spatial models. In the finite-dimensional setting, the asymptotic behavior of the maxima of Gaussian random vectors was first investigated by Huesler and Reiss \cite{HR89}. Recently,  Kabluchko, Schlather and de Haan   \cite{A1} and Kabluchko \cite{A7} consider the functional setting and prove convergence of the maxima of i.i.d. centered Gaussian processes under suitable conditions on their covariance structure.
The family of  Brown-Resnick processes has quickly become a popular model for spatial extremes due to its simple characterization by a negative definite functions. Our purpose here is to revisit their results and apply the present framework: we put the emphasis on regular variations and point processes. In the sequel, we follow the approach by Kabluchko, see Theorems 2 and 6 in \cite{A7}.

We suppose that $(T,d)$ is a compact metric space satisfying the entropy condition
\[
 \int_{0}^{1} (\log N ( \varepsilon ))^{1/2} d \varepsilon < \infty, 
\]
where $N ( \varepsilon ) $ is the smallest number of balls of radius needed to cover $T$. \\
Let  $Z_{n}=\{Z_n(t);\ t\in T\}$, $n\geq 1$, be a sequence of continuous zero-mean  unit-variance Gaussian processes with covariance function 
\[
r_{n} ( t_{1} , t_{2}) = \mathbb E [ Z_{n} (t_{1} ) Z_{n} ( t_{2}) ], \quad t_{1}, t_{2} \in T.
\]
Define the scaling sequence
\begin{equation}\label{eq:scaling}
b_{n}=(2 \log n)^{1/2} - (2 \log n )^{-1/2} ( ( 1/2 )\log \log n + \log (2 \sqrt{\pi})),\quad n\geq 1
\end{equation}
and the rescaled process
\[
  Y_{n}(t) = b_{n} ( Z_{n}(t) - b_{n} ), \quad t \in T. 
\]
We consider the log-normal process 
\[
  X_{n} (t) = \exp Y_{n}( t), \quad t \in T.
\]

\begin{theo}\label{theokabluch} Fix $t_{0} \in T$ and suppose that:
\vspace{-0.4cm}
\begin{enumerate}[ i )]
 \item\label{condi} Uniformly in $t_{1}, t_{2} \in T$ : 
\begin{equation}\label{eqlem1}
\Gamma ( t_{1}, t_{2} ) = \lim_{n \rightarrow \infty} 4\log n ( 1 - r_{n} ( t_{1} , t_{2}) )\  \in [0, \infty ).
\end{equation} 
\item  For all  $t_{1}, t_{2} \in T$, there exists  $C > 0$ such that 
\begin{equation}\label{eqlem2} 
\sup\limits_{ n \geqslant 1 }\log n  ( 1 - r_{n} ( t_{1} , t_{2} ) )  \leqslant C d( t_{1} , t_{2}).
\end{equation} 
\end{enumerate}
Then, 
\[
\lim_{n \rightarrow \infty} n \mathbb P [ X_{n} \in \cdot\mbox{ }]  \stackrel{w ^{\sharp} }{\longrightarrow} \nu ( \cdot) \quad \mbox{ in }  M_{b}^{\sharp} (\overline{ \mathbb C}_0^+ )  , 
\]
where 
\[
\nu (A) =\int _{ 0} ^{ \infty }  \mathbb{P}  [ w e^{W(\cdot) -\frac{1}{2} \Gamma(t_{0},\cdot) } \in A ] w^{-2} dw,\quad   A \subseteq \mathbb C_0^+\ \mbox{Borel\ set} ,
\]
and $\lbrace W(t)  , t \in T \rbrace$ a centered Gaussian process  such that $W(t_{0})=0$ and with incremental variance $\Gamma$.
\end{theo}

Following \cite{A7,A1}, we define $Y_{n}^{w}$ the process $Y_{n}$ conditioned by $Y_{n} (t_{0}) = w$ and we note  $\mu_{n}^{w} (t)=\mathbb E [Y_{n}^{w}(t)]$. 
We will need the following  two Lemmas.

\begin{lemme}\label{lemfinn1}
Under conditions i) and ii) of Theorem \ref{theokabluch}, the family $\{ Y_{n}^{w} ;\ n \geq 1 \}$ is tight in $\bbC(T)$ for all fixed $\omega\in\bbR$, as well as the family $\{Y_{n}^{w} - \mu_{n}^{w};\  w \in \mathbb R ,  n \geq 1\}$.
\end{lemme}

Lemma \ref{lemfinn1} follows from the proof of  Theorem 6 in \cite{A7} (see also the proof of Theorem 17 in \cite {A1}). Details are omitted here.  

\begin{lemme}\label{lem5}
Under condition i) of Theorem \ref{theokabluch}, the following convergence in the sense of finite dimensional distributions holds: as $n\to\infty$,
 \[
 \{  Y_{n}^{w} (t) ,\ t \in T \}   \stackrel{f.d.d.}{\longrightarrow}  \{ w + W(t) - \frac{1}{2}\Gamma ( t , t_{0} ) ,\ t \in T\} 
 \]
with $W$ defined in Theorem \ref{theokabluch}
\end{lemme}

\begin{preuve}{of Lemma \ref{lem5}  }\\
Standard computations for Gaussian processes entails that $Y_n^w$ is a Gaussian process with mean and covariance 
\begin{equation}\label{eq:mean}
 \mu_{n}^{w}(t) = w   r_{n} (t,t_{0})  + b_{n}^{2}( r_{n} (t,t_{0}) -1),
\end{equation}
\begin{equation}\label{eq:var}
 \mbox{Cov} [ Y_{n}^{w}(t_{1}),  Y_{n}^{w}(t_{2})  ] = b_{n}^{2}  \left( r_{n}   ( t_{1} , t_{2}) -r_{n}   ( t_{1} , t_{0})   r_{n} ( t_{2} , t_{0}) \right) .
\end{equation}
Under assumption i), it holds
$$ \lim_{ n \rightarrow \infty} \mu_{n}^{w}(t) = w-\frac{1}{2}\Gamma(t,t_{0})$$
and 
$$ \lim_{ n \rightarrow \infty} \mbox{Cov} [ Y_{n}^{w}(t_{1}), Y_{n}^{w} (t_{2})  ] =\frac{1}{2} ( \Gamma (t_{1}, t_{0}) + \Gamma (t_{2}, t_{0}) - \Gamma (t_{1}, t_{2}) ).$$
This implies that $\mbox{Var}(Y_{n}^{w}(t_{2})- Y_{n}^{w} (t_{1}))\to\Gamma(t_1,t_2)$ and proves the Lemma.
\end{preuve}

\begin{preuve}{  of Theorem \ref{theokabluch} }\\
Our proof relies on the following criterion for $\sharp$-convergence in $\overline\bbC_0^+$  by  Hult et Lindskog. According to Theorem 4.4 of \cite{A5}, it is enough to prove that:
\begin{enumerate}[i)]
\vspace{-0.4cm}
\item for every $r>0$, $\sup \limits_{n\geq 1} n \mathbb P [ \sup_{t \in K} X_{n}(t) > r] < \infty$;
\item for every $\varepsilon>0$, $ \lim \limits _{ \delta \rightarrow 0 } \sup\limits_{n\geq 1} n \mathbb P [  \omega_{ X_{n}} (\delta) \geqslant \varepsilon] = 0$,\\
 where $\omega_{f}(\delta)=\sup\{|f(t_1)-f(t_2)|;\ t_1,t_2\in T, d(t_1,t_2)\leq\delta\}$ denotes the modulus of continuity of $f\in\bbC(T)$;
\item $n\bbP[X_n \in\,\cdot\,]\to \nu$ in the sense of finite-dimensional convergence.
\end{enumerate}
These three points are proven below. At several places, we will use that the scaling sequence $b_n$ defined by Equation \eqref{eq:scaling} satisfies
$\sqrt{2\pi} b_{n} e^{b_{n}^{2}/2} \sim  n$ as $n\to\infty$. As a consequence, the sequence $\dfrac{n}{\sqrt{2\pi} b_{n} e^{b_{n}^{2}/2}}$ is bounded by some constant $M>0$.

\begin{itemize}
\vspace{-0.4cm}
\item[$\bullet$] Proof of i):\\
Let $r>0$ and $\tilde{r} = \ln r$. We have
\begin{eqnarray*}
  n \mathbb P [  \sup \limits_{t \in T} X_{n}(t) > r] 
&=&   n \mathbb P [ \sup_{t \in T} Y_{n}(t) >  \tilde r] \\ 
&=&\frac{n}{\sqrt{2\pi} b_{n} e^{b_{n}^{2}/2} } \int _{ \bbsR } e^{- w - w^{2} / 2 b_{n}^{2}}  \mathbb{P}  \lbrack \sup \limits_{t \in T} Y_{n}^{w} (t) > \tilde{r}  \rbrack  dw   \\ 
& \leq& M \int _{ \bbsR } e^{- w }  \mathbb{P}  \lbrack \sup \limits_{t \in T} Y_{n}^{w} (t) > \tilde{r} \rbrack  dw,  
\end{eqnarray*} 
so that
\begin{equation}\label{point1}
\sup_{n\geq 1}  n \mathbb P [  \sup \limits_{t \in T} X_{n}(t) > r] 
\leq M \int _{ \bbsR } e^{- w }  \sup_{n\geq 1} \mathbb{P}  \lbrack \sup \limits_{t \in T} Y_{n}^{w} (t) > \tilde{r} \rbrack  dw . 
\end{equation} 
First, we have  
\begin{equation}\label{point1.1}
\int _{0}^{+ \infty} e^{- w }  \sup_{n\geq 1}\mathbb{P}  \lbrack \sup_{t \in T} Y_{n}^{w} (t) > \tilde{r}  \rbrack  dw\leqslant  \int _{0}^{+ \infty} e^{- w }   dw  < \infty. 
\end{equation} 
By tightness of the family $ \{Y_{n}^{w} - \mu_{n}^{w};\  w \in \mathbb R ,  n \geq 1\}$ (see Lemme \ref{lemfinn1}), there exists  $c_{1}> 0$ such that for all $n\geq 1$ and $w\in\bbR$, 
$$ \mathbb{P}  \lbrack \sup_{t \in T} ( Y_{n}^{w} (t) - \mu_{n} ^{w} ( t) ) > c_{1}\rbrack < \frac{1}{2}. $$
Furthermore, Equations \eqref{eq:mean}-\eqref{eq:var} and assumption ii) of Theorem \ref{theokabluch} together imply that there are some  $c_{2},c_3>0$ and $n_0\geq 1$ such that for all $n\geq n_0$ and all $w<0$,
$$ \sup_{t \in T}  \mu_{n} ^{w} ( t) \leqslant \frac{1}{2}w + c_{2} , \mbox{ and } \sup_{t \in T}  \mbox{Var}[ Y_{n} ^{w} ( t) ]\leqslant c_{3}^{2}.$$
Applying Borell's inequality (see Theorem D.1 in \cite{P96}), we obtain
$$ \mathbb{P}  \lbrack \sup_{t \in T} Y_{n}^{w} (t) > \tilde{r} \rbrack  < 2 \psi \Big(  -\frac{\tilde{r}   -w/2 -c_{1} -c_{2}}{c_{3}}\big)\leq 2e^{-(\frac{\tilde{r}   -w/2 -c_{1} -c_{2}}{c_{3}})^{2}},\quad  w<0 ,$$
where $\psi$ is the tail of the standard Gaussian distribution. 
Consequently,  
\begin{equation}\label{point1.2}
\int ^{0}_{- \infty} e^{- w }\sup_{n\geq n_0}  \mathbb{P}  [ \sup_{t \in T} Y_{n}^{w} (t) > \tilde{r} ] dw < 2 \int ^{0}_{- \infty} e^{- w}e^{-(\frac{\tilde{r}   -w/2 -c_{1} -c_{2}}{c_{3}})^{2}}   dw<\infty.
\end{equation}
Equations  \eqref{point1.1} and \eqref{point1.2} together imply 
\begin{equation}\label{point1.3}
 \int_{\bbsR} e^{- w}\sup_{n\geq 1}\mathbb{P}  [ \sup_{t \in T} Y_{n}^{w} (t) > \tilde{r} ] dw<\infty.
\end{equation}
In view of Equation \eqref{point1}, this proves i).
 
\item[$\bullet$] Proof of ii):\\
For $\delta>0$, we have
\begin{eqnarray*}
  \sup\limits_{n\geq 1} n \mathbb P [  \omega_{X_{n}} (\delta) \geqslant \varepsilon] &=& \sup\limits_{n\geq 1}  \dfrac{n}{\sqrt{2\pi} b_{n} e^{b_{n}^{2}/2} } \int _{ \bbsR } e^{- w - \frac{w^{2} }{2 b_{n}^{2}}}  \mathbb{P}  [ \omega_{e^{Y_{n}^{w} }} (\delta) \geqslant \varepsilon ]  dw\\
 & \leq & M\int _{ \bbsR } e^{- w } \sup\limits_{n\geq 1}  \mathbb{P}  [ \omega_{e^{Y_{n}^{w} }} (\delta) \geqslant \varepsilon ]  dw. 
\end{eqnarray*}
We will apply Lebesgue dominated convergence Theorem to prove that this has limit zero as $\delta\to 0$. 
Lemma \ref{lemfinn1} implies that for every fixed $w$, the sequence $ \{  e^{Y_{n}^{w} };\ n \geqslant 1 \}$ is tight so that
\[
\lim \limits _{ \delta \rightarrow 0 } e^{-w} \sup\limits_{n\geq 1}  \mathbb{P}  [ \omega_{e^{Y_{n}^{w} }} (\delta) \geqslant \varepsilon ] = 0.
\]
It remains to prove a suitable domination condition: since $\delta\mapsto \mathbb{P}  [ \omega_{e^{Y_{n}^{w} }} (\delta) \geqslant \varepsilon ]$ is a nondecreasing function, it is enough to prove that for some $\delta_0>0$
\[
\int _{ \bbsR } e^{- w } \sup\limits_{n\geq 1}  \mathbb{P}  [ \omega_{e^{Y_{n}^{w} }} (\delta_0) \geqslant \varepsilon ]  dw<\infty.
\]
Since $\omega_{e^{Y_{n}^{w} }} (\delta_0)\leq \exp [\sup_{t\in T} Y_n^w(t)]$, we have 
\[
e^{- w } \sup\limits_{n\geq 1} \mathbb{P}  [ \omega_{e^{Y_{n}^{w} }} (\delta_{0}) \geqslant \varepsilon ]  \leqslant e^{- w } \sup\limits_{n\geq 1} \mathbb{P}  [ \sup_{t\in T} Y_{n}^{w}(t)  \geqslant  \ln \varepsilon ]. 
\]
Equation \eqref{point1.3} with $\tilde r=\ln \varepsilon$ provides the required domintation condition.

\item[$\bullet$] Proof of iii):\\
We have to prove that for all set $A\subseteq\overline C_0^+$ of the form 
\[
A=\{f\in\bbC_0^+; \exists j\in [\![1,k]\!],\ f(t_j)> x_k\},\quad k\geq 1,\ t_1,\ldots,t_k\in T,\ x_1,\ldots,x_k\in (0,+\infty),
\]
it holds
\[
 \lim_{n\to\infty} n\bbP[X_n\in A] =\nu(A).
\]
Letting  $y_{j}= \ln x_{j} $,  we have
\begin{eqnarray*}
 n \mathbb{P}[X_n\in A]&=& n\bbP  \lbrack \exists j \in [\![1,k]\!],\ X_{n} (t_{j}) > x_{j}   \rbrack \\
 & =& n \mathbb{P}  \lbrack \exists j \in [\![1,k]\!],\ Y_{n} (t_{j}) > y_{j}   \rbrack   \\ 
 &=&\dfrac{n}{\sqrt{2\pi} b_{n} e^{b_{n}^{2}/2} } \int _{\bbsR } e^{- w - \frac{w^{2} }{2 b_{n}^{2}}}  \mathbb{P}  \lbrack \exists j \in [\![1,k]\!],\ Y_{n}^{w} (t_{j}) > y_{j}  \rbrack  dw   \\ 
 &\sim &\int _{\bbsR}e^{- w - \frac{w^{2} }{2 b_{n}^{2}}}  \mathbb{P}  \lbrack \exists j \in [\![1,k]\!],\ Y_{n}^{w} (t_{j}) > y_{j}  \rbrack  dw.   
\end{eqnarray*}
We will apply once again Lebesgue's dominated convergence Theorem.  Lemma \ref{lem5} entails 
\[
\mathbb{P}  \lbrack \exists j \in [\![1,k]\!],\ Y_{n}^{w} (t_{j}) > y_{j}  \rbrack \to \mathbb{P}  \lbrack \exists j \in [\![1,k]\!],\ w+W(t_j)-\frac{1}{2}\Gamma(t_j,t_0) > y_{j} \rbrack,
\]
so that we expect 
\begin{eqnarray*}
 n \mathbb{P}[X_n\in A]&\to& \int _{\bbsR}e^{- w }  \mathbb{P}  \lbrack \exists j \in [\![1,k]\!],\ w+W(t_j)-\frac{1}{2}\Gamma(t_j,t_0) > y_{j} \rbrack  dw\\
&=&\int _{0}^\infty  \mathbb{P}  \lbrack \exists j \in [\![1,k]\!],\ we^{W(t_j)-\frac{1}{2}\Gamma(t_j,t_0)} > x_{j} \rbrack w^{-2} dw \\
&=&\nu(A).  
\end{eqnarray*}
We are left to prove a suitable domination condition: for all $n\geq 1$, we have 
\[
 e^{- w - \frac{w^{2} }{2 b_{n}^{2}}}  \mathbb{P}  \lbrack \exists j \in [\![1,k]\!],\ Y_{n}^{w} (t_{j}) > y_{j}  \rbrack 
\leq e^{-w}\sup_{n\geq 1} \mathbb{P}  \lbrack \sup_{t\in T} Y_{n}^{w} (t) > \tilde r  \rbrack
\]
with $\tilde r=\min_{1\leq j\leq k} y_j$. Equation \eqref{point1.3}  yields the required domination condition.
\end{itemize}
\end{preuve}



\bibliographystyle{plain}
\bibliography{Biblio}

\end{document}